\documentclass[10pt]{article}
\usepackage[dvips]{graphics}
\usepackage{epsfig,graphicx}
\usepackage{amsfonts}
\usepackage{amsmath}

\input xy
\xyoption{all}

\newcommand{\lar}{\longrightarrow}
\newcommand{\kF}{\mathcal F}

\newcommand{\kG}{\mathcal G}

\newcommand{\kB}{\mathcal B}
\newcommand{\kM}{\mathcal M}
\newcommand{\kO}{\mathcal O}

\newcommand{\kJ}{\mathcal J}
\newcommand{\kE}{\mathcal E}
\newcommand{\kT}{\mathcal T}
\newcommand{\kN}{\mathcal N}
\newcommand{\kS}{\mathcal S}
\newcommand{\kk}{k}
\newcommand{\kA}{\mathcal A}
\newcommand{\PP}{\mathbf P}
\newcommand{\ZZ}{\mathbb Z}

\newcommand{\FF}{\mathbf F}
\newcommand{\tens}{\otimes}
\newcommand{\Hom}{\mathop{\rm Hom}}
\newcommand{\Coh}{\mathop{\rm Coh}}
\renewcommand{\mod}{\mathop{\rm mod}}
\newcommand{\Ext}{\mathop{\rm Ext}\nolimits^1}
\newcommand{\TC}{\mathop{\rm TC}}
\newtheorem{theorem}{Theorem}[section]
\newtheorem{corollary}[theorem]{Corollary}
\newtheorem{remark}[theorem]{Remark}
\newtheorem{example}[theorem]{Example}
\newtheorem{lemma}[theorem]{Lemma}
\newtheorem{statement}[theorem]{Statement}
\newtheorem{definition}[theorem]{Definition}
\newtheorem{conjecture}[theorem]{Conjecture}

\begin{document}

\title
{
Coherent Sheaves on Singular Projective Curves with Nodal Singularities
 }

\author
{
Igor Burban \and Yurij Drozd}
\date{}

\maketitle

\begin{abstract}
 We give the full answer to the question: on which curves the category
of coherent sheaves $\Coh_{X}$ is tame. The answer is: these are just the curves
from the list of Drozd-Greuel (see ~\cite{vb}). Moreover, in this case the 
derived category
$D^{-}(\Coh_{X})$ is also tame.  We give an explicit description of the objects 
of this category as well as of the categories $D^{b}(\Coh_{X})$, $\Coh_{X}$.
Among the coherent sheaves we describe the vector bundles, torsion-free
sheaves, mixed sheaves and skyscraper sheaves.
\end{abstract}

\section{\bf Introduction}

Let $X$ be a projective curve over $k=\bar{\kk}$. For any two coherent sheaves
$\kF$ and $\kG$ we have 
$$
\dim_{\kk}(\Hom(\kF,\kG)) < \infty.
$$
This implies that in the category of coherent sheaves the generalized 
Krull-Schmidt theorem holds:
$$
  \kF \cong \bigoplus\limits_{i=1}^{s} \kF_{i}^{m_{i}},
$$ 
where $\kF_{i}$ are indecomposable and $m_i, \kF_{i}$ are uniquely defined.

\noindent
\underline{Our aim}: to describe all indecomposable coherent sheaves on $X$. 
Consider first a smooth case. For any coherent sheaf $\kF$ we have 
$$
0 \longrightarrow T(\kF) \longrightarrow \kF \longrightarrow \kF/T(\kF) 
\longrightarrow 0,
$$
where $T(\kF)$ is the torsion part of $\kF$ (skyscraper sheaf) and 
$\kF/T(\kF) $ the torsion-free quotient of $\kF$. But $\kF/T(\kF)$ is even
locally free, since our curve is smooth.
The local-global spectral sequence 
$H^{p}({\mathcal Ext}^{q}(\ , \ )) 
\Longrightarrow \mathop{\rm Ext}\nolimits^{p+q}(\ , \ )$
implies 
$$
\hspace{-1.1cm}
0 \lar H^{1}({\mathcal Hom}(TF(\kF/T(\kF),T(\kF)) \lar
\Ext(\kF/T(\kF),T(\kF)) \lar
H^{0}({\mathcal Ext}^{1}(\kF/T(\kF),T(\kF)) \lar 0.
$$

So, $\Ext(\kF/T(\kF),T(\kF))=0$ and $\kF\cong T(\kF))\oplus \kF/T(\kF)$.
Also, in a smooth case, indecomposable objects of $\Coh_{X}$ are 
\begin{itemize}
\item scyscraper sheaves ${\mathcal O_{x}/m_{x}^{n}}$
\item indecomposable vector bundles.
\end{itemize}

\noindent What is known about the classification of indecomposable vector 
bundles on smooth projective curves?
\begin{enumerate}
\item Let $X = \PP^{1}_{k}$. Then indecomposable vector bundles are just 
the line bundles $\kO_{\PP^{1}}(n), n \in \ZZ$ (see ~\cite{Groth}).
\item Let $X$ be an elliptic curve. The indecomposable vector bundles are 
described by two discrete parameters $r,d$: rank and degree and one continuous
(point of the curve $X$)  ~\cite{At}.

\item It is well-known that with the growth of the genus of the curve $g(X)$ 
the moduli spaces of 
vector bundles become bigger and bigger. For the smooth curves of genus 
$g \ge 2$ it was shown (Drozd/Greuel ~\cite{vb}, Scharlau(1992)) 
that the problem
of classification of vector bundles is wild. ``Wild'' means

\begin{enumerate}
\item ``geometrically'': we have $n$-parameter families of indecomposable vector bundles for arbitrary large $n$;
\item "algebraically": for every finite-dimensional $\kk$-algebra $\Lambda$
there is an exact functor $(\Lambda-\mod) \lar VB_{X}$ mapping non-isomorphic
objects to non-isomorphic ones and indecomposable into indecomposable.

\end{enumerate}

\end{enumerate}

Moreover, Drozd and Greuel have proved the following trichotomy (see ~\cite{vb}):

\begin{enumerate}
\item $VB_{X}$ is finite (indecomposable objects are described by discrete parameters)
if $X$ is a configuration of projective lines of the type $A_n$.

\begin{figure}[ht]
\hspace{5.5cm}
\includegraphics[height=2.5cm,width=1.7cm,angle=-90]{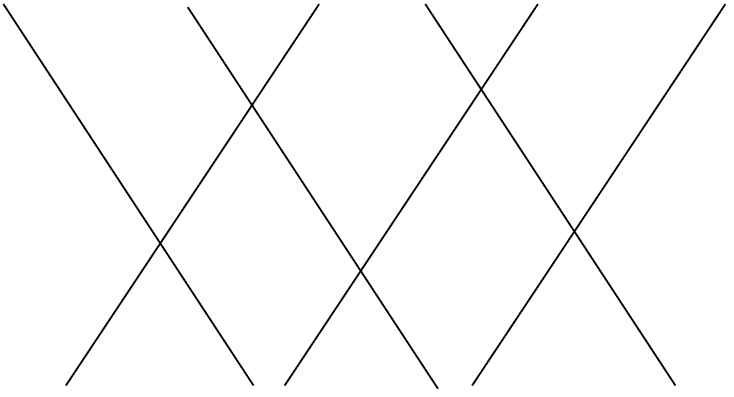}
\end{figure}

\item $VB_{X}$ is tame (intuitively this means that indecomposable objects are 
pa\-ra\-metri\-zed by 
1 continuous parameter and several discrete parameters. For the rigorous definition see also ~\cite{vb}) if
\begin{enumerate}
\item $X$ is an elliptic curve

\begin{figure}[ht]
\hspace{5.5cm}
\includegraphics[height=1.5cm,width=2cm,angle=-90]{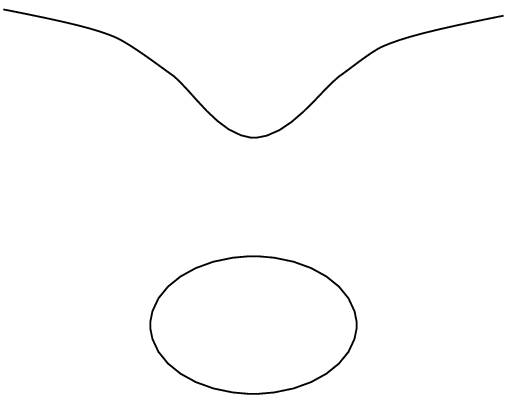}
\end{figure}

\clearpage
\item $X$ is a rational curve with one simple node
\begin{figure}[ht]
\hspace{5.5cm}
\includegraphics[height=1.5cm,width=2cm,angle=-90]{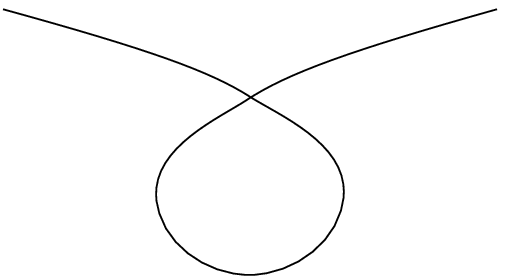}
\end{figure}

\item $X$ is a configuration of projective lines of type $\tilde{A_n}$

\begin{figure}[ht]
\hspace{4.5cm}
\includegraphics[height=3cm,width=2cm,angle=-90]{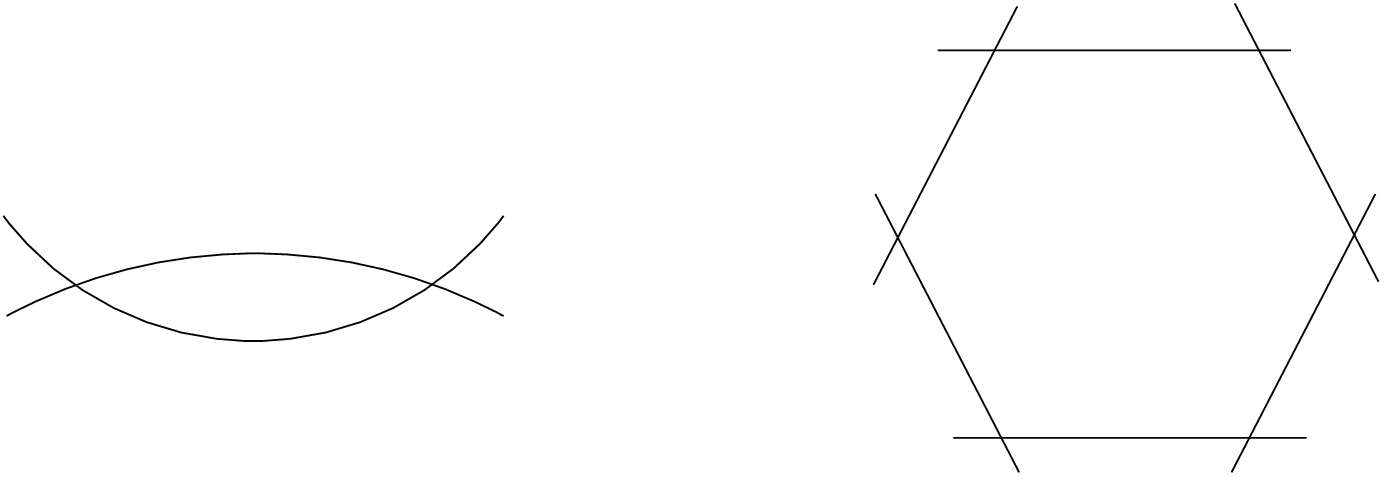}
\end{figure}

\end{enumerate}

\item Wild, otherwise.
\end{enumerate}

Drozd and Greuel also gave an explicit description of vector bundles and
torsion-free sheaves. We want to describe the coherent sheaves in the cases
$1-2$. The problem is that the sequence
$$
0 \longrightarrow T(\kF) \longrightarrow \kF \longrightarrow \kF/T(\kF) 
\longrightarrow 0
$$
does not necessarily split. ($\kF/T(\kF)$ could be torsion-free but not locally 
free.) We can descibe all skyscraper sheaves (see ~\cite{Roiter1},~\cite{Roiter2},~\cite{GP}) and 
torsion-free sheaves. But there is no common technology of reconstructing the objects possibly being 
at the middle. We need a new idea and this idea is to 
apply the technique of the derived categories. There is a full and faithful 
functor
$\Coh_{X} \hookrightarrow D^{-}(\Coh_{X})$ sending a coherent sheaf $\kF$ into its 
locally free resolution $\kF_{.}$. We shall describe indecomposable objects of the bigger
category  $D^{-}(\Coh_{X})$ and among them the complexes with zero higher
homologies (which correspond to coherent sheaves).

{\it Aknowledgement}. We would like to thank the working group of algebraic geometry of the Univesity
of Kaiserslautern and especially Prof. G.-M.~Greuel for the hospitality 
and helpful discussions. We are also grateful to P.~Bitsch for the 
typographical correction of this article. 

\section{\bf Main construction}

Let $X$ be a projective curve of the type as in the list
of Drozd-Greuel 
(singular points are simple double points or transversal intersections), 
$\tilde{X} \stackrel{\pi}\lar X$ its normalization, $\tilde{\kO}=
\pi_{*}(\kO_{\tilde X})$ (so $\tilde{\kO}_{x} = \bar{\kO}_{X,x} $ --- the
integral closure of $\kO_{X,x}$), $\kJ = Ann_{\kO}(\tilde{\kO}/\kO)$ the
conductor (so the support of $\kO/\kJ$ is precisely the set of singular points 
$Sing(X)$).
Since the morphism $\pi$ is affine, we can identify 
$\kO_{\tilde X}$-modules and $\tilde \kO $-modules. 

\begin{definition} Consider the following category of triples of 
complexes $\TC_{X}$ (for convenience, all objects of the derived categories 
below are supposed to be the complexes of locally-free modules)

\begin{enumerate}
\item Objects are the triples 
$(\tilde {\mathcal F}_{.}, {\mathcal M}_{.}, i)$, where

$\tilde{\mathcal F}_{.} \in D^{-}(\Coh_{\tilde{\kO}})$, 

${\mathcal M}_{.} \in D^{-}(\Coh_{\kO/{\mathcal J}})$, 

$i : {\mathcal M}_{.} \lar
 \tilde{\mathcal F}_{.} \tens_{\tilde \kO} \tilde{\kO}/{\mathcal J}$ 
( $\tens_{\tilde \kO} \tilde{\kO}/{\kJ}$ denotes the derived functor of the
tensor product) a morphism in
 $ D^{-}(\Coh_{\kO/{\mathcal J}}) $, such that 

$i\tens id :
{\mathcal M}_{.} \tens_{\kO} \tilde{\kO}/{\mathcal J} \lar
\tilde{\mathcal F}_{.} \tens_{\tilde \kO} \tilde{\kO}/{\mathcal J}$
is an isomorphism in $D^{-}(\Coh_{\tilde{\kO}/\kJ})$ (we implicitly use here that $gl.dim(\Coh_{\kO/\kJ})=0$. 
The last means that the morphism $i\tens id$ 
is correctly defined). We want to stress (it is the only exception to the agreement above)
that $i\tens id$ is {\it
not the derived functor} of the tensor product but just the tensor product.

\item Morphisms
 $                              
(\tilde{\mathcal F_{1}}_{.}, {\mathcal M_{1}}_{.}, i_{1}) \lar
(\tilde{\mathcal F_{2}}_{.}, {\mathcal M_{1}}_{.}, i_{2})
 $
are  pairs $(\Phi, \varphi)$, 
$\tilde{\mathcal F_{1}}_{.} \stackrel{\Phi}\lar \tilde{\mathcal F_{2}}_{.},
 {\mathcal M_{1}}_{.} \stackrel{\phi}\lar
 {\mathcal M_{2}}_{.}
$, such that 

\begin{tabular}{p{2.5cm}c}
 &
\xymatrix
{
\tilde{\mathcal F_{1}}_{.} \ar[r] \ar[d]^{\Phi} & 
\tilde{\mathcal F_{1}}_{.} \tens_{\tilde \kO} \tilde{\kO}/{\kJ} 
\ar[d]^{
\Phi\tens id} &
{\mathcal M_{1}}_{.} \ar[l]_{\qquad i_1} \ar[d]^{\phi} \\                                                                            
\tilde{\mathcal F_{2}}_{.} \ar[r] & 
\tilde{\mathcal F_{2}}_{.} \tens_{\tilde \kO} \tilde{\kO}/{\kJ} &
{\mathcal M_{2}}_{.} \ar[l]^{\qquad i_2} 
}
\end{tabular}

is commutative (more precisely, we want the right square to be commutative in
$D^{-}(\Coh_{\kO/\kJ})$).
\end{enumerate}
\end{definition}

\begin{theorem}
The functor 
$$
D^{-}({\rm Coh}_{X}) \stackrel{\mathbf F}\lar {\rm TC}_{X}
$$
$
{
\mathcal
 F^{.} \lar
({\mathcal F}_{.} \tens_{\kO} \tilde{\kO},
{\mathcal F}_{.} \tens_{\kO} \kO/{\mathcal J},
i : {\mathcal F}_{.} \tens_{\kO} \kO/{\mathcal J}   \lar {\mathcal F}_{.} \tens_{\kO} {\tilde \kO}/ {\mathcal J})
}
$
fulfils the following properties:

\begin{enumerate}
\item ${\mathbf F}$ is dense (i.e., every triple $(\tilde{\kF}_{.}, \kM_{.}, i)$ is 
isomorphic to some
${\mathbf F}(\kF_{.})$). 
\item $\FF(\kF_{.}) \cong \FF(\kG_{.}) \Longleftrightarrow 
\kF_{.} \cong \kG_{.}$.
\item ${\mathbf F}$ is full.
\end{enumerate} 
\end{theorem}

\begin{remark} $\FF$ is not faithful. So it is not an equivalence.
\end{remark}

\noindent
{\it Proof.}
The main ingredient of the proof is: having a triple 
$(\tilde{\kF}_{.}, \kM_{.}, i)$ how can we reconstruct $\kF_{.}$?

\noindent
The exact sequence 
$$0\lar \kJ\tilde{\kF}_{.} \lar \tilde{\kF}_{.} \lar
\tilde{\kF}_{.}\tens_{\tilde{\kO}} \tilde{\kO}/\kJ \lar 0$$ 
in $\Coh_{X}$ gives a distinguished triangle 
$$
\kJ\tilde{\kF}_{.} \lar \tilde{\kF}_{.} \lar
\tilde{\kF}_{.}\tens_{\tilde{\kO}} \tilde{\kO}/\kJ 
\lar \kJ\tilde{\kF}_{.}[-1]
$$ 
in $D^{-}(\Coh_{X})$.
The properties of the triangulated categories imply the morphism of triangles

\begin{tabular}{p{1.4cm}c}
&
\xymatrix
{
{\mathcal J}\tilde{\mathcal F}_{.} \ar[r] &
\tilde{\mathcal F}_{.} \ar[r] &
\tilde{\mathcal F}_{.} \tens_{\tilde \kO} \tilde{\kO}/{\mathcal J} \ar[r]&
{\mathcal J}\tilde{\mathcal F}_{.}[-1] \\
{\mathcal J}\tilde{\mathcal F}_{.} \ar[r] \ar[u]^{id} &
{\mathcal F}_{.} \ar[r] \ar[u]_{\Phi} &
{\mathcal M}_{.} \ar[r] \ar[u]_{i} &
{\mathcal J}\tilde{\mathcal F}_{.}[-1]. \ar[u]^{id}\\
}
\end{tabular}

In other words $\kF_{.} = cone(\kM_{.} \lar \kJ \kF_{.}[-1])[1]$. Taking a cone
is not a functorial operation. It gives an intuitive explanation, why
functor ${\bf F}$ is not an equivalence.

The properties of the triangulated categories imply immediately that
the constructed map 
$$
{\bf G} : Ob({\rm TC}_{X}) \lar Ob(D^{-}({\rm Coh}_{X}))
$$
sends isomorphic objects into isomorphic ones and satisfies ${\bf GF}(\kF_{.}) \cong \kF_{.}$. 
Now we have to show, that
${\bf FG}(\tilde{\kF}_{.}, \kM_{.}, i) \cong (\tilde{\kF}_{.}, \kM_{.}, i)$.
Suppose we have found a representative of $\kM_{.}$ (in the corresponding class
of the derived category) such that 
$i\tens id $ is an isomorphism {\it of complexes}. Consider the pull-back diagram
in the abelian category $Com(\Coh_{X})$:

\begin{tabular}{p{1.6cm}c}
&
\xymatrix
{0 \ar[r] & \kJ \tilde{\kF}_{.} \ar[r] \ar[d]^{id} & \kF_{.} \ar[r]^{\Psi} \ar[d]^{\Phi} &
 \kM_{.} \ar[r] \ar[d]^{i} & 0 \\
0 \ar[r] & \kJ \tilde{\kF}_{.} \ar[r]  & \tilde{\kF}_{.} \ar[r]^{\pi}  &
 \tilde{\kF}_{.}/ \kJ \tilde{\kF}_{.} \ar[r] & 0.
}
\end{tabular}

\noindent
We are going to show that
\begin{enumerate}
\item $\kF_{.}$ is a complex of locally-free $\kO$-modules;
\item $(\Phi\tens id, \Psi\tens id) : (\kF_{.}\tens_{\kO}\tilde{\kO}, \kF_{.}\tens_{\kO} \kO/\kJ,
\kF_{.}\tens_{\kO} \kO/\kJ \lar \kF_{.}\tens_{\kO}\tilde{\kO/\kJ}) \lar (\tilde{\kF}_{.}, \kM_{.}, i)$
is an isomorphism in the category of triples. 
\end{enumerate}
The first condition should be checked in each stalk of every component of the complex $\kF_{.}$. We have to show
that $\forall x \in X, n\ge 0 \ (\kF_{n})_{x}$ is a locally-free $\kO_{X,x}$-module. In a regular point it is 
obvious. Let $x$ be a singular point of a curve. Denote
$
A=\kO_{X,x}, \tilde{A}= \tilde{\kO}_{X,x}, J=\kJ_{x}, F=(\kF_{n})_{x}, 
\tilde{F}= (\tilde{\kF_{n}})_{x}, M=(\kM_{n})_{x}.
$
From the projectivity of $A$ follows a commutative diagram:

\begin{tabular}{p{3cm}c}
&
\xymatrix
{
& & & A^r \ar[d] \ar[ddl]_{\psi} & \\
& & & (A/J)^r \ar[d]^{iso} & \\
0 \ar[r] &J\tilde F \ar[r] \ar[d]^{id} & F \ar[r]^{\Psi} \ar[d]^{\Phi} &
M \ar[r] \ar[d]^{i} & 0 \\
0 \ar[r] &J\tilde{F} \ar[r]  & \tilde{F} \ar[r]  &
\tilde{F}/J\tilde{F} \ar[r] & 0. \\
}
\end{tabular}

Since $\Phi$ is injective, $F$ is a torsion-free module. 
Tensor the first row of the exact sequence   with $\tens_{A} \tilde{A}/J$ and the second with 
$\tens_{\tilde A} \tilde{A}/J$.
Since $i\tens id $ is an isomorphism,  we obtain that $\Phi\tens id$ is an isomorphism modulo the radical. 
By Nakayama's lemma,
$\Phi\tens id$ is an epimorphism. By the same reason $\psi$ is also an epimorphism.
But $rank(F)=rank(\tilde{F})=r$. $\psi $ is an epimorphism of torsion-free modules of the same rank. So, $\psi$ 
is an isomorphism. The same holds for $\Phi\tens id$.
We have shown, moreover, that $\Phi\tens id$ is an isomorphism of complexes. 
The same we can say about $\Psi\tens id$. Now let us prove the statement 
we have used.

\begin{statement} Let $\tilde{\kF}_{.}$ be a complex of locally-free $\tilde{\kO}$-modules, $i: \kM_{.} \lar \tilde{\kF}_{.}/\kJ$
a morphism in $D^{-}(\Coh_{\kO/\kJ})$ such that $i\tens id$ is an isomorphism in $D^{-}(\Coh_{\tilde{\kO}/\kJ})$. Then there 
is a complex $\kM'_{.}$ and a homotopy isomorphism $f : \kM'_{.}\lar \kM_{.}$ such that for $i'=if$  the map
$i'\tens id : M'\tens_{\kO/\kJ} \tilde{\kO}/\kJ \lar \tilde{\kF}_{.}/\kJ$ is an isomorphism of complexes. 
\end{statement}

\noindent
{\it Proof.}
Let us mention the following easy 
\begin{lemma}
Let  $\kA$ be an abelian category of a homological dimension 0, $X_{.}$
and $Y_{.}$ two complexes,  $f_{.} : X_{.} \lar Y_{.}$ a morphism of complexes.
Then the following holds:
\begin{enumerate}
\item Let $B_{n}(X)=Im(d_{n+1}(X_{.}))$. Then $X_{n} = B_{n}(X) \oplus H_{n}(X) 
\oplus B_{n-1}(X)$. 
\item Let $d_{n} : B_{n}(X) \oplus H_{n}(X) \oplus B_{n-1}(X) \lar
B_{n-1}(X) \oplus H_{n-1}(X) \oplus B_{n-2}(X)$ be given by 
$d_n(b_n,h_n,b_{n-1})=(b_{n-1},0,0)$. It defines a complex, quasi-isomorphic to  
$X_{.}$. 
\item In such a notation the map $f_{.}$ looks as follows:
 $$
f_n=
\left(
\begin{array}{ccc}
f_n|_{B_n(X)} & 0 & 0 \\
0 & H_{n}(f) & 0 \\
0 & 0& f_{n-1}|_{B_{n-1}(X)}
\end{array}
\right).
$$
\item Moreover, the morphism of complexes 
 $$
\tilde{f_n}=
\left(
\begin{array}{ccc}
0 & 0 & 0 \\
0 & H_{n}(f) & 0 \\
0 & 0& 0
\end{array}
\right)
$$
is homotopic to $f_{.}$. We'll call it also a canonical form of a morphism $f$.
\end{enumerate}
\end{lemma}

Suppose $X$ is a curve with simple nodes.
 Consider complex $\tilde{\kF}_{.}/\kJ$ from
$D^{-}(\Coh_{\tilde{\kO}/\kJ})$. We see 
that all $(\tilde{\kF}_{n}/\kJ)_{x}$ ($x$ is a singular point of $X$) considered as $k$-vector spaces are 
even-dimensional. Moreover, 
all  homologies $H_{n}((\tilde{\kF}_{.}/\kJ)_{x})$ are even-dimensional, too.
So, all the boundaries $B_{n}((\tilde{\kF}_{.}/\kJ)_{x})$
are also even-dimensional. The proof is an easy linear algebra now, so we shall
leave it.
We have shown that our functor ${\bf F}$ is dense and is a bijection on the 
iso-classes of 
indecomposable objects. Let us show that {\bf F}
is full.

Let $(\Phi, \phi) : (\tilde{\kF_{1}}_{.}, {\kM_{1}}_{.}, i_1) \lar 
(\tilde{\kF_{2}}_{.}, {\kM_{2}}_{.}, i_2)$ be a morphism in 
$TC_{X}$, where $i_1, i_2$ induce an isomorphism of complexes after tensoring,
$\Phi$  {\it is a morphism of complexes}. If both $\phi$ and $\bar{\Phi}$ are 
in the canonical form, then the following 
 diagram  

\begin{tabular}{p{2.7cm}c}
&
\xymatrix
{ {\mathcal M_{1}}_{.} \ar[r]^{\phi} \ar[d]^{i_1} & {\mathcal M_{2}}_{.} 
\ar[d]^{i_2} \\
  \tilde{\mathcal F_{1}}_{.}/ \kJ\tilde{\kF_{1}}_{.} \ar[r]^{\bar{\Phi}} &
  \tilde{\kF_{2}}_{.}/ \kJ\tilde{\kF_{2}}_{.} \\
}
\end{tabular}

\noindent
is commutative 
{\it in the category of complexes}. The properties of pull-back imply the 
existence of a 
morphism of complexes 
$ {\kF_{1}}_{.} \lar {\kF_{2}}_{.} $ such that  
$$
\begin{tabular}{p{1.3cm}c}
&
\xymatrix@!0
{ & {\kF_{2}}_{.} \ar[rr]\ar'[d][dd] & & {\kM_{2}}_{.} \ar[dd] \\
 {\kF_{1}}_{.} \ar[ur] \ar[rr] \ar[dd]& & {\kM_{1}}_{.} \ar[ur]
\ar[dd] \\
& {\tilde{\kF_{2}}}_{.} \ar[rr] & 
& \tilde{\kF_{2}}_{.}/ \kJ\tilde{\kF_{2}}_{.} \\
\tilde{\kF_{1}}_{.} \ar[rr] \ar[ur] & 
& \tilde{\kF_{1}}_{.}/ \kJ\tilde{\kF_{1}}_{.} \ar[ur]
}
\end{tabular} 
$$
is commutative. If $\Phi : \tilde{\kF_{1}}_{.} \lar \tilde{\kF_{1}}_{.}$ 
is a
quasi-isomorphism, then 
$\bar{\Phi} : \tilde{\kF_{1}}_{.}/ \kJ \tilde{\kF_{1}}_{.}  \lar 
\tilde{\kF_{2}}_{.}/ \kJ\tilde{\kF_{2}}_{.}$
is a quasi-isomorphism, too
(the tensor product is a functor). The axioms of triangulated categories imply 
that  $\kJ\tilde{\kF}_{.}  \lar \kJ\tilde{\kG}_{.} $  is also a quasi-isomorphism. 
Hence, the induced map ${\kF_{1}}_{.} \lar 
{\kF_{2}}_{.}$  
is  also a quasi-isomorphism. It implies that the functor ${\bf F}$ is full.
Nevertheless it is not faithful. It can be seen as follows: let 
 
$$ {\mathcal E}_{.}[-1]= \dots \lar 0 \lar \underbrace{\kE}_{1} \lar 0 \dots 
$$
 and 
 $$\kF_{.} = \dots \lar 0 \lar \underbrace{\kF}_{0} \lar 0 \dots $$
 be two complexes 
($\kF$ and $\kG$ are vector bundles). 
$\Hom(\kF_{.}, \kE_{.}[1])= \Ext(\kF,\kE)$. The map 
$$ 
\hspace{-0.2cm}
\mathop{\rm Hom}\nolimits_{\scriptsize D^{-}({\rm Coh}_{X})}(\kF_{.}, \kE_{.}[1])=
\Ext_{\kO_{X}}(\kF,\kE) \lar  \Ext_{\tilde{\kO_{X}}}(\tilde{\kF},\tilde{\kE})
=\mathop{\rm Hom}\nolimits_{{\rm TC}_{X}}({\bf F}(\kF), {\bf F}(\kG))$$ is not a monomorphism: 
 $\Ext(\kO,\kF)=H^{1}(\kF)$, but
$H^{1}(X,\kF) \ne H^{1}(\tilde{X},\tilde{\kF})$. 
So our functor is not faithful.
It proves our theorem.

\begin{corollary}
In the notation as above let $\kE$ and $\kG$ be two vector bundles on $X$.
Then  the canonical map 
$$
\Ext_{\kO_{X}}(\kF,\kE) \lar  \Ext_{\tilde{\kO}_{X}}(\tilde{\kF},\tilde{\kE})
$$
is surjective. 
\end{corollary}

\section{\bf Coherent sheaves on a rational curve with one node}
Let us consider first the case of the rational curve with one simple node.
Suppose its equation is $zy^{2} - x^{3} - x^2 z = 0$. Then its normalization
is 
$\tilde{X} = {\bf P}^1$. Suppose that the preimages of a singular point are 
$(0:1)=0$ and $(1:0)=\infty$.

\begin{figure}[ht]
\hspace{4.5cm}
\includegraphics[height=3.5cm,width=2cm,angle=-90]{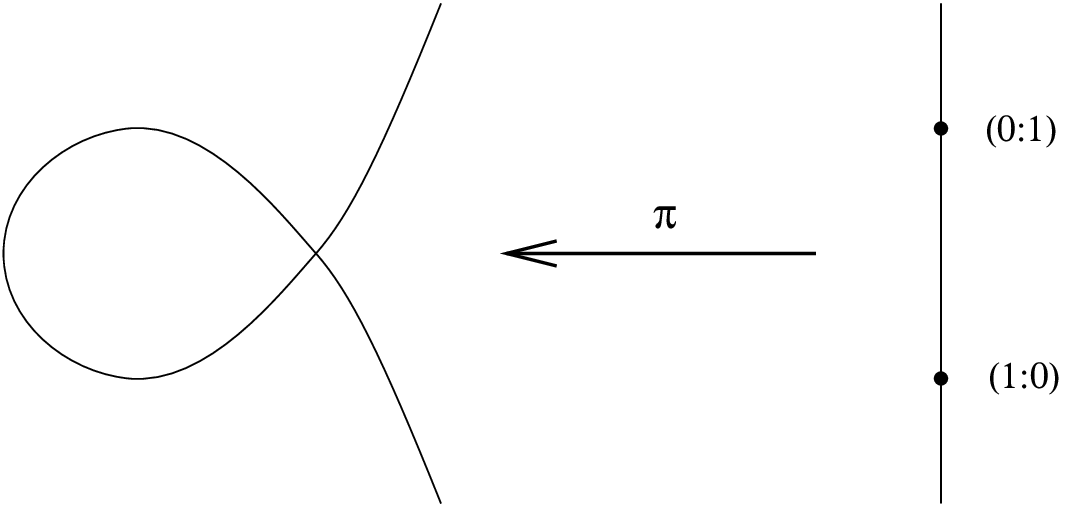}
\end{figure}

What does the result of the previous section mean? A complex $\kF_{.}$ from
the derived category $D^{-}(\Coh_{X})$ as a data structure is uniquely defined 
by some triple $(\tilde{\kF}_{.}, \kM_{.}, i)$. 

What is $\tilde{\kF}_{.}$?
The category $\Coh_{\bf P^1}$ has the global dimension $1$. It means
(Dold, 1960) that indecomposable objects of $D^{-}(\Coh_{\bf P^1})$ are
$$
\kE_{n}[r] : \dots \lar 0 \lar \underbrace{\kO_{\bf P^1}(n)}_{r} \lar 0 \lar 0 \lar \dots,
$$
$$
\kT_{kx}[s] : \dots\lar 0 \lar \kO_{\bf P^1}(-kx) \hookrightarrow 
\underbrace{\kO_{\bf P^1}}_{s} \lar 0 \lar \dots .
$$
A complex $\tilde{\kF}_{.}$ is just a direct sum
$$
   \tilde{\kF}_{.} \cong \bigoplus(\kE_{n}[r] \oplus \kT_{kx}[s]).
$$

Now let us explain what is $\kM_{.}$ and what is $i$. $\kO/\kJ$ is a 
skyscraper sheaf $k_{0}$ (with the stalk $k$ at the singular point), 
$\tilde{\kO}/\kJ = (k\times k)_{0}$. 
It means that the category $\Coh_{\kO/\kJ}$ is semi-simple. So, 
$\kM_{.} \cong (H_{.}(\kM_{.}),0))$ and we get, moreover, a commutative diagram

\begin{tabular}{p{3cm}c}
&
\xymatrix
{ \kM_{.} \ar[r]^{i} \ar[d] & \tilde{\kF}_{.}/\kJ) \ar[d] \\
   H_{.}(\kM_{.}) \ar[r]^{H_{.}(i)} & H_{.}(\tilde{\kF}_{.}/\kJ)
}
\end{tabular} 

The map 
$H_{k}(i) :  
H_{k}(\kM_{.}) \lar H_{k}(\tilde{\kF}_{.}/\kJ)$ is simply a map of two 
vector spaces. But  $H_{k}(\tilde{\kF}_{.}/\kJ)$ is also a  $k\times k$-module.
This implies that  $H_{k}(i)$ is given by two matrices 
$H_{i}(0)$ and $H_{i}(\infty)$
 (intuitively  --- one corresponds to the point 
$0$, another to  $\infty$). Moreover, both of these matrices have the same
size and are nondegenerated. It follows from the condition $3$ of 
the category of triples:
$ H_{k}(\kM_{.})\tens_{\kO/\kJ} \tilde{\kO}/\kJ 
\xrightarrow{ H_{k}(i)\tens id} H_{k}(\tilde{\kF_{.}}/ \kJ)
$ 
is an isomorphism.
 Choose a set of local parameters of $\tilde{\kF}$ in each component of the complex $\kF_{.}$. They induce some basis
in $H_{.}(\tilde{\kF})$.  Choose some basis in $H_{k}(\kM_{.})$. With respect
to such choice the map 
 $H_{.}(i)$ is given by a collection  of matrices 

\begin{figure}[ht]
\hspace{-2.4cm}
\includegraphics[height=17cm,width=4.5cm,angle=-90]{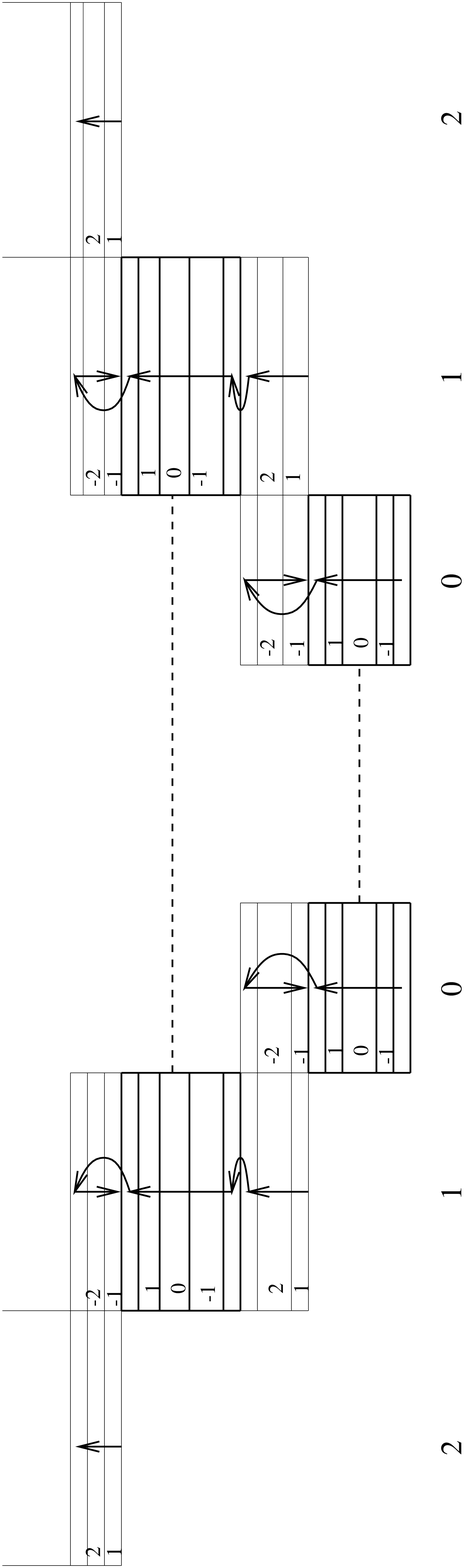}
\end{figure}



There are two types of blocks --- those, which came from vector bundles and those
which came from skyscraper sheaves.  The blocks are numbered by integers and
natural numbers, respectively. This numbering defines some ``weights'' of the 
blocks of vertical matrices (the partial order is shown in the picture).
 Blocks, corresponding to the same skyscraper
or vector bundle, are called conjugated. Conjugated blocks have the same number 
of rows.  Indeed, all the blocks of $H_{i}(0)$ and $H_{i}(\infty)$, but finitely many have the zero size. 
But if one of the conjugated blocks is nonempty, then the other one is 
nonempty, too.

 Now we should answer on the following
question: which triples $(\tilde{\kF}_{.},\kM_{.},H_{.}(i))$ correspond to the 
isomorphic complexes $\kF_{.}$? Surely, we have to consider the automorphisms
of $\tilde{\kF}_{.}$ and look at what they induce in homologies.  As a result 
we get the following matrix problem:

\begin{enumerate}
\item we can do any simultaneous elementary transformations of the columns of 
the matrices
 $H_{i}(0)$ and $H_{i}(\infty)$;
\item we can do any simultaneous transformations of rows inside conjugated blocks;
\item we can add a scalar multiple of any row from the block 
with lower weight to any row of a block of the higher weight (inside of the big matrix, of course). These transformations can be proceeded independently
inside of $H_{i}(0)$ and $H_{i}(\infty)$ (see the next section for  more details).
\end{enumerate}






These types of matrix problems are well-known in the representation theory.
First they appeared in the work of Nazarova-Roiter (~\cite{Roiter1}) about the 
classification of $k[[x,y]]/(xy)$-modules. They are called,
sometimes, Gelfand problems in honour of I.~M.~Gelfand, who 
formulated a 
conjecture (at the International Congress of Mathematics in 
Nice (1970))
about the structure of Harish-Chandra modules at the singular point of 
$SL_{2}({\mathbb R})$ ~\cite{Gelfand}. 
This problem was reduced to some matrix problem of such a type ~\cite{Roiter3}.

The strict categorical formulation and then a solution of this type of problem was done 
by Nazarova-Roiter, Drozd, 
Bondarenko ~\cite{Roiter3},~\cite{bimproblems},~\cite{mp1}. It means that these matrices correspond to
the objects of some category and these objects are isomorphic if and only if
 one matrix can be transformed into another one by the above set of transformations. 
From this point of view, it is enough to describe the indecomposable objects.

Let us recall  a combinatoric of the answer in this case. There are two types
of indecomposable objects: bands and strings ~\cite{mp}. 
We shall give the definitions in the next section. Here we only want to stress, 
that band object depends on one continuous parameter and several discrete
parameters.
String object depands only on discrete parameters. Let us  consider some 
examples.

\begin{example}
The following data (band)  define a simple vector bundle of rank 2 on 
$X$:
normalization $\tilde{\kO}\oplus \tilde{\kO(n)}, n\ne 0$, matrices:
\begin{figure}[ht]
\hspace{5cm}
\includegraphics[height=2.5cm,width=1cm,angle=-90]{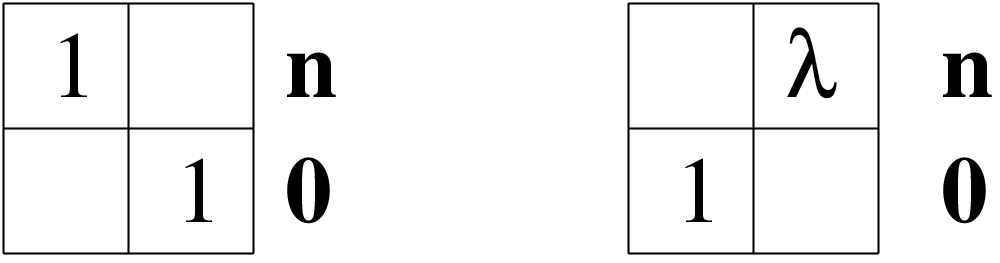}
\end{figure}
\end{example}

\begin{example}
The following data (band) define a mixed sheaf (sheaf which is neither torsion 
nor torsion-free): normalization is 
$$
\kE_{0}[0]^{2n}\oplus \kE_{-3}[0]^{n} \oplus
\kT_{2 0}[0]^{n} \oplus \kT_{5 0}[0]^{n} \oplus
\kT_{1 \infty}[0]^{n} \oplus \kT_{4 \infty}[0]^{n},
$$
matrices are  
\begin{figure}[ht]
\hspace{0.7cm}
\includegraphics[height=10cm,width=2.8cm,angle=-90]{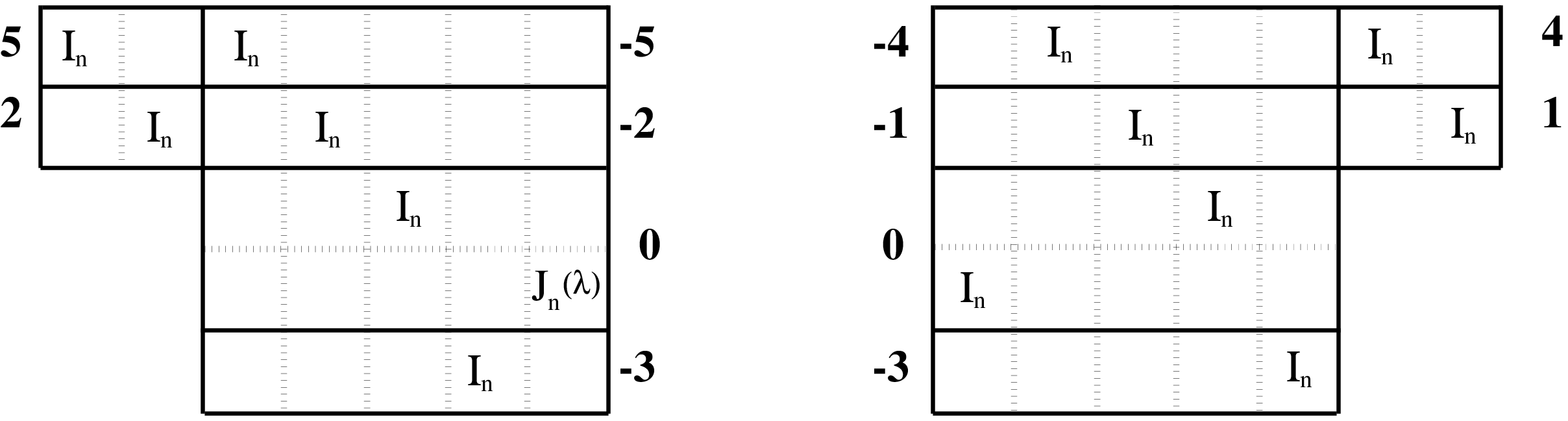}
\end{figure}
\end{example}

\begin{example}
The following data (string) define the skyscraper sheaf $k_{0}$:
normalization is 
$$\bigoplus\limits_{i=0}^{\infty}(\kT_{1 0}[-i] \oplus 
\kT_{1 \infty}[-i]),$$
matrices:
\clearpage 
\begin{figure}[ht]
\hspace{2.5cm}
\includegraphics[height=6.5cm,width=2.8cm,angle=-90]{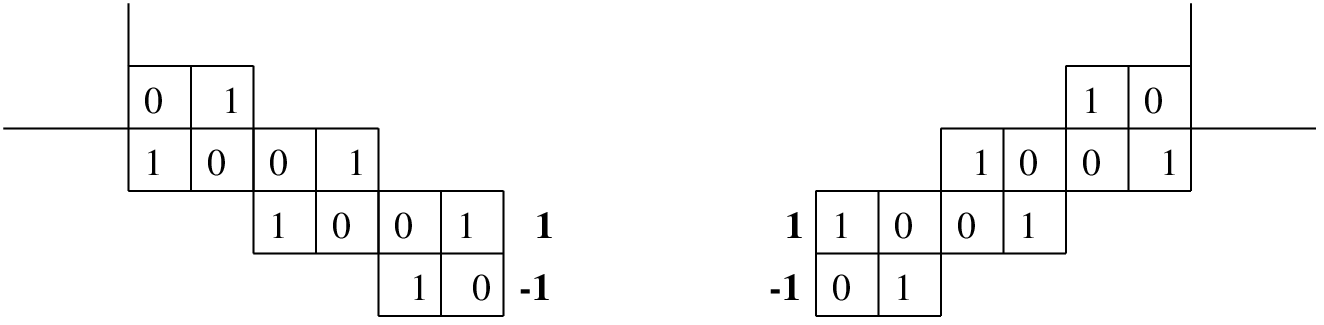}
\end{figure}
\end{example}

\begin{example}
The following data  (band) defines some object of the bounded
category $D^{b}(Coh_{X})$ which is not a coherent sheaf:
normalization is 
$$
\hspace{-0.7cm}
(\kE_{-2}[-2])^n \oplus (\kT_{2 0}[-1])^n \oplus (\kT_{1 \infty}[-1])^n
\oplus (\kE_{0}[-1])^n \oplus (\kT_{1 0})^{2n} \oplus (\kT_{2 \infty})^n 
\oplus (\kT_{4 \infty})^n \oplus (\kE_{-1})^n 
\oplus (\kE_{1})^n
$$
matrices:
\begin{figure}[ht]
\includegraphics[height=12cm,width=3.3cm,angle=-90]{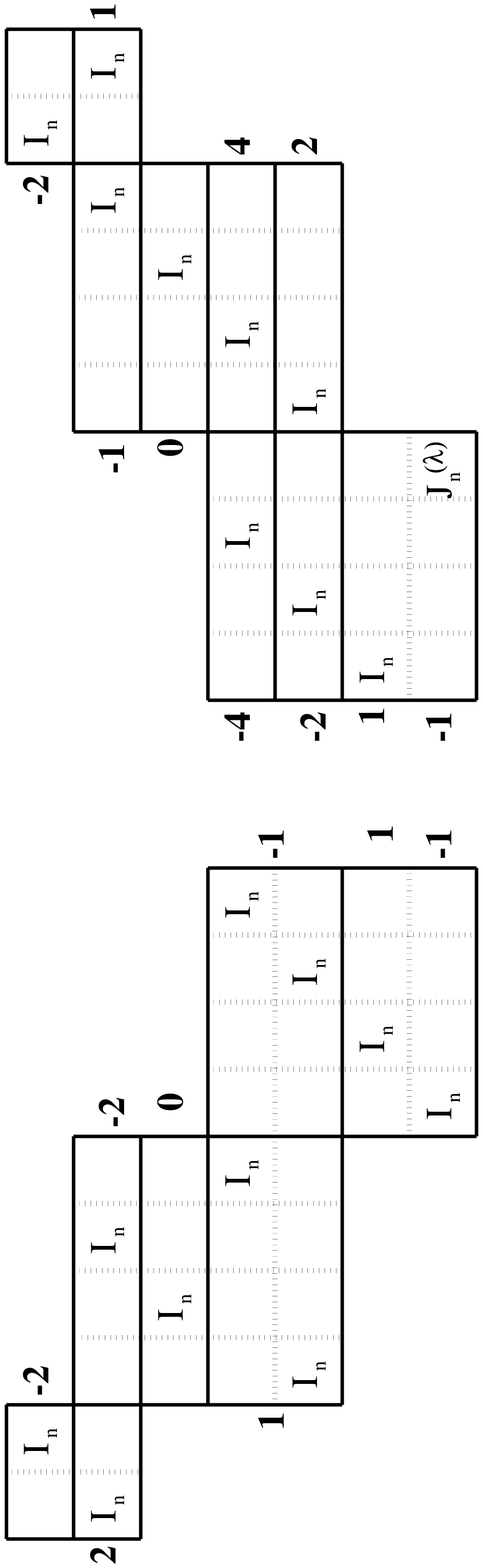}
\end{figure}
\end{example}

\begin{example}
The following data  (string) defines some object of the 
category $D^{-}(Coh_{X})$ which is not an object of  the bounded
derived category  $D^{b}(Coh_{X})$:
normalization is
$$
\dots \oplus \kT_{1 0}[-2] \oplus \kT_{1 \infty}[-2] \oplus \kE_{2}[-2] \oplus
\kT_{3 0}[-1] \oplus \kT_{1 \infty}[-1] \oplus \kE_{-1}[-1] \oplus 
\kT_{2 0} \oplus \kT_{3 \infty} \oplus \kE_{0}
$$
matrices:
\begin{figure}[ht]
\includegraphics[height=12cm,width=3.3cm,angle=-90]{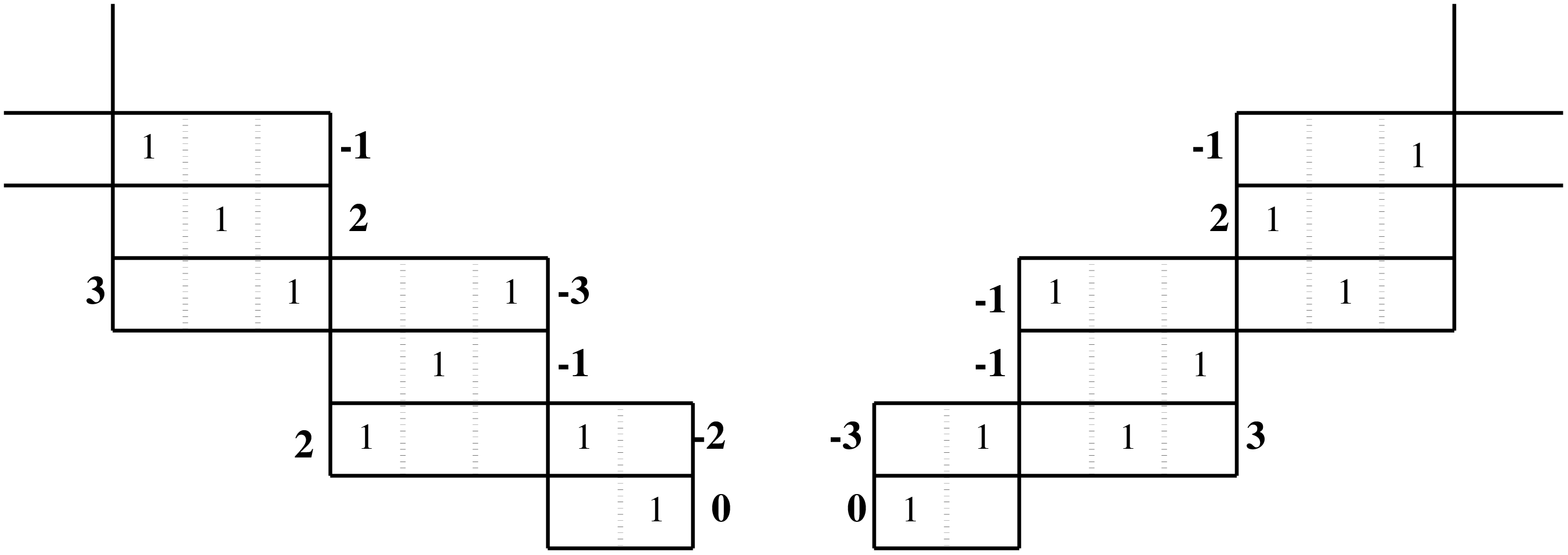}
\end{figure}
\end{example}

Let us now consider a general case.

\section
{\bf Reduction to the matrix problem}

Suppose $X$ is either a rational curve with one node or a configuration of
projective lines of the type $A_n$, $\tilde{A_n}$, $\tilde{X} \stackrel{\pi}
\lar X $ its normalization. Then $\tilde{X}$ is just a disjoint union of 
projective lines,
$\tilde{X} = \coprod\limits_{i=1}^{n} X_{i}$. Let $\pi_{i}$ be the restriction 
of $\pi$ on $X_{i}$, $\tilde{\kO}_{i}= (\pi_{i})_{*}(\pi^{*}(\kO_{X}))$,
$\kJ_{i}=\pi_{*}(\kJ|_{X_{i}})$. 
The category $\Coh_{\tilde X}$ is 
equivalent to $\prod\limits_{i=1}^{n} \Coh_{X_{i}}$ (the same holds
for the derived categories),
$\tilde{\kO}-\mod$ is equivalent  to the direct product of the categories
$(\tilde{\kO}_{1}-\mod) \times (\tilde{\kO}_{2}-\mod) \times \dots 
\times (\tilde{\kO}_{n}-\mod)$.
$\tilde{\kO}/\kJ-\mod $ is equivalent  to 
$(\tilde{\kO}_{1}/\kJ_{1}-\mod) \times (\tilde{\kO}_{2}/\kJ_{2}-\mod) 
\times \dots \times (\tilde{\kO}_{n}/\kJ_{n}-\mod)$.
Choose the local coordinates on each of the lines $X_{i}$ in such a way
that the preimages of the singular points are either $(0 : 1)$ or
$(1 : 0)$. We can interpret coherent $\tilde{\kO}_{i}$-modules just as 
 $\Coh_{\bf P^1}$.
$\kJ_{i}$ will be identified with  the ideal sheaf of two points $(0 : 1)$ and $(1 : 0)$, or, if only one point lies on $X_{i}$, just with the ideal sheaf of $(0 : 1)$. 
Both $\kO/\kJ$ and $\tilde{\kO}/\kJ$ are the skyscraper sheaves with support in the singular points of $X$. Suppose $X$ is the configuration of
projective lines of the type $\tilde{A}_{n}$. The canonical morphism 
$\kO/\kJ \lar \tilde{\kO}/\kJ$ is then the diagonal  morphism
$$
k\times k \times \dots \times k \lar (k\times k) \times (k\times k)\times \dots
\times (k\times k).
$$

Let $(\tilde{\kF}_{.}, \kM_{.}, i)$ be some triple.  $\tilde{\kF}_{.} \cong
\tilde{\kF_{1}}_{.}\oplus \tilde{ \kF_{2}}_{.}\oplus \dots \oplus 
\tilde{ \kF_{n}}_{.}$,
where $\tilde{\kF_{i}}_{.} \in D^{-}(\Coh_{X_{i}})= D^{-}(\Coh_{\bf P^1})$.

How do we get a matrix problem in this case? 
$\Coh_{\mathbf P^1}$ has a homological dimension $1$, which implies that every
indecomposable object of $D^{-}(\Coh_{\mathbf P^1})$ is isomorphic to some object
of the type
$\dots \lar 0 \lar \underbrace{\kF}_{i} \lar 0\lar \dots$, where 
$\kF \in Ob(\Coh_{\mathbf P^1})$  is indecomposable. Indecomposable 
objects of  
 $\Coh_{\bf P^1}$
are known: line bundles 
$\kO_{\bf P^1}(n)$ and skyscraper sheaves $\kO_{{\bf P^1},x}/m_{x}^n$. 

The skyscraper sheaf 
$\kO_{{\bf P^1},x}/m_{x}^n$ has a locally free resolution 
$$
0\lar \kO_{\bf P^1}(-n) \lar \kO_{\bf P^1} 
\lar \kO_{{\bf P^1},x}/m_{x}^n \lar 0,
$$
which means that in the derived category 
$ (\kO_{\bf P^1}(-n) \lar \kO_{\bf P^1}) \cong \kO_{x}/m_{x}^n $ holds.

Choose the local bases in each component of the complex $\tilde{\kF}$.
The map
$H_{k}(i) : H_{k}(\kM_{.}) \lar H_{k}(\tilde{\kF}_{.}/\kJ)$ is given
by $n+1$ matrices, corresponding to the singular points of $X$. Each of 
these matrices  itself consists of two nondegenerated components of the same
size.

The morphisms in the derived category $D^{-}(\Coh_{\bf P^1})$ and their 
images after tensoring
with $\tilde{\kO}/\kJ$ and taking the homology, is very easy to compute in this 
case, so
we just leave this computation.

The question is, which transformations can we do with the matrices defining the 
homology?
From the definition of the category of triples follows that we have to consider
the automorphisms $\Phi : \tilde{\kF}_{.} \lar \tilde{\kF}_{.}$ and
$\phi :  \kM_{.}\lar \kM_{.}$ which make the following diagram
$$
\begin{tabular}{p{1.3cm}c}
&
\xymatrix
{
H_{.}(\kM_{.}) \ar[r]^{H_{.}(i)} \ar[d]^{H_{.}(\phi)} &
H_{.}(\tilde{\kF}_{.}/\kJ) \ar[d]^{H_{.}(\Phi\tens id)} \\
H_{.}(\kM_{.}) \ar[r]^{H_{.}(i)}&
H_{.}(\tilde{\kF}_{.}/\kJ)  \\
}
\end{tabular}
$$
commutative.

With every curve from the list of Drozd-Greuel we can associate some partially
ordered set.
Let  $\omega_{-1}< \omega_{0} < \omega_{1}$ be three cardinal numbers
(it means that $n\omega_{-1}<m\omega_{0} < k \omega_{1} \ \forall
n,m,k \in \ZZ$).
The algorithm is the following now: consider the set of pairs $(L,a)$,
where $L$ is a component of $\tilde{X}$, $a \in L$ some preimage of the 
singular point.
To such a pair we associate the following partially ordered sets:
$F_{(L,a)}(k)$, $E_{(L,a)}(k)$ $(k \in {\mathbb N})$. $F_{(L,a)}(k)$ 
consist from one element.
$E_{L.a}(0)$ has two types of elements
$E_{(L,a)}(0, m\omega_{-1})$, $m \in \ZZ_{-}$ and 
$E_{(L,a)}(0, n\omega_{0}), n\in \ZZ_{+}$.
 $E_{(L,a)}(i) (i \ge 1)$ has three types of elements: 
$E_{(L,a)}(i, m\omega_{-1}) (m \in \ZZ_{-})$,
$E_{(L,a)}(i, n\omega_{0}) (n \in \ZZ)$, $E_{(L,a)}(i, k\omega_{1}) (k \in \ZZ_{+})$. 
Each set $E_{(L,a)}(i)$ has a natural partial order.

Consider the union of all points 

$$
{\bf E\bigcup F} = (\bigcup\limits_{L,a,i} E_{L,a}(i)) \bigcup (\bigcup\limits_{L,a,i} F_{L,a}(i)).
$$

In this set let us introduce the following equivalence relation:

\begin{enumerate}
\item $E_{(L,a)}(i, -n\omega_{-1}) \sim E_{(L,a)}(i+1, n\omega_{1})$, $i \ge 0$;
\item $E_{(L,a)}(i, m\omega_{0}) \sim E_{(L',a)}(i, m\omega_{0})$ $i \ge 0$;
\item $F_{(L,a)}(i) \sim F_{(L,a')}(i)$ $i \ge 0$.
\end{enumerate}

\begin{example}
Rational curve with one node 
\begin{figure}[ht]
\hspace{4cm}
\includegraphics[height=3.5cm,width=4cm,angle=-90]{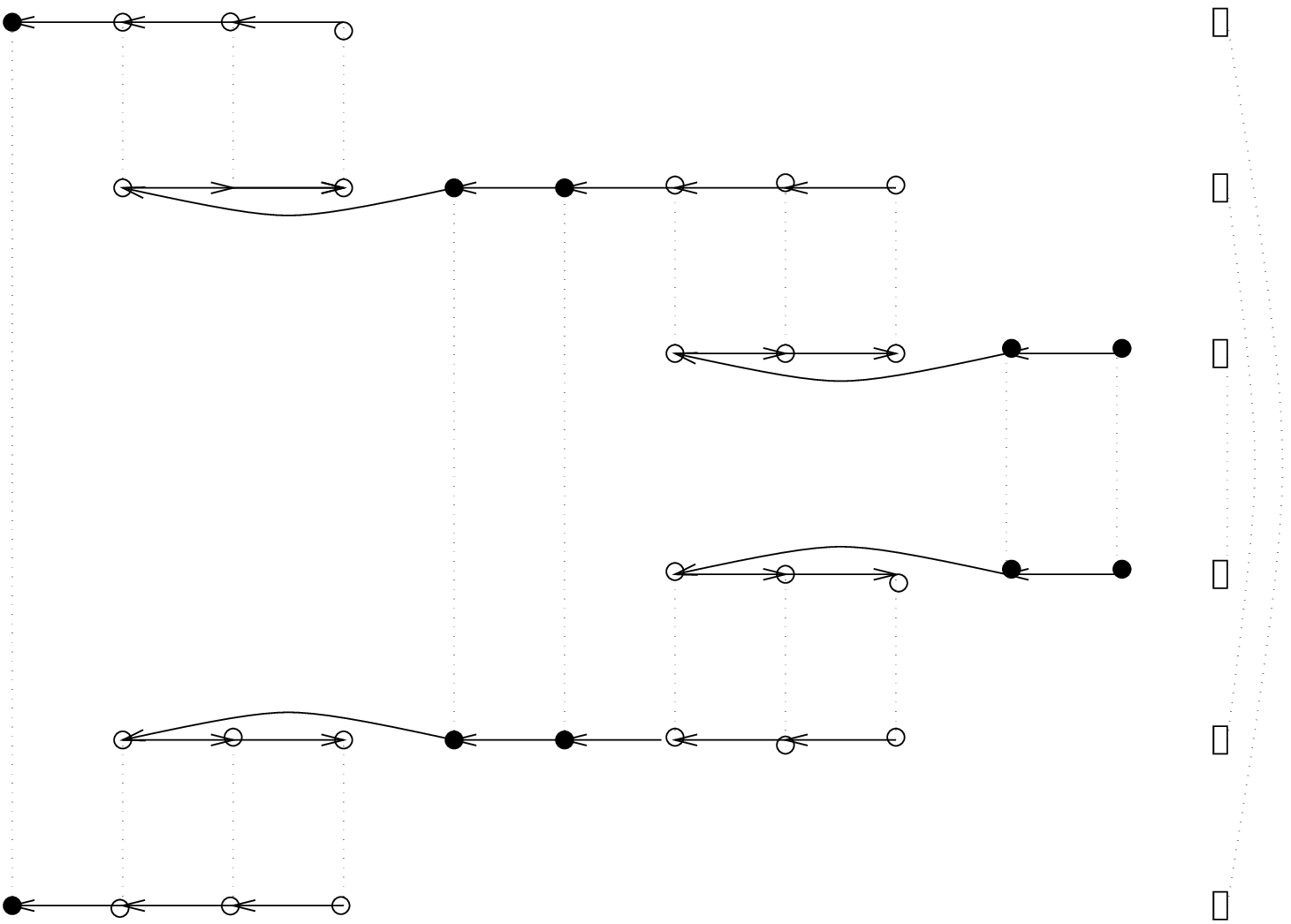}
\end{figure}
\end{example}

\clearpage
\begin{example}
Transversal intersection of two lines at one point ($A_{1}$-case)
\begin{figure}[ht]
\hspace{4cm}
\includegraphics[height=3.5cm,width=4cm,angle=-90]{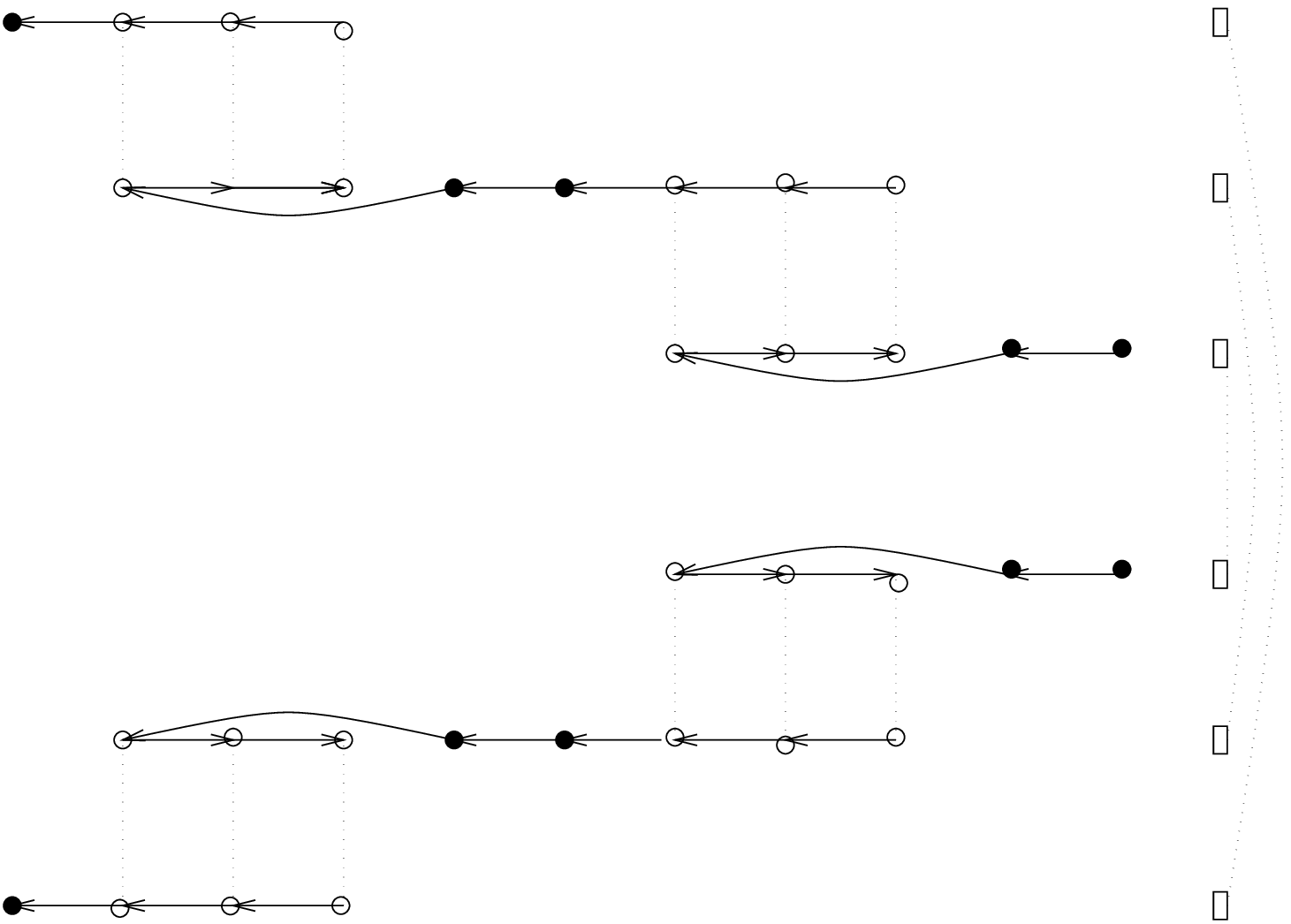}
\end{figure}
\end{example}

\begin{example}
Transversal intersection of two lines at 
two points ($\tilde{A}_{1}$-case) 
\begin{figure}[ht]
\hspace{4cm}
\includegraphics[height=3.5cm,width=6cm,angle=-90]{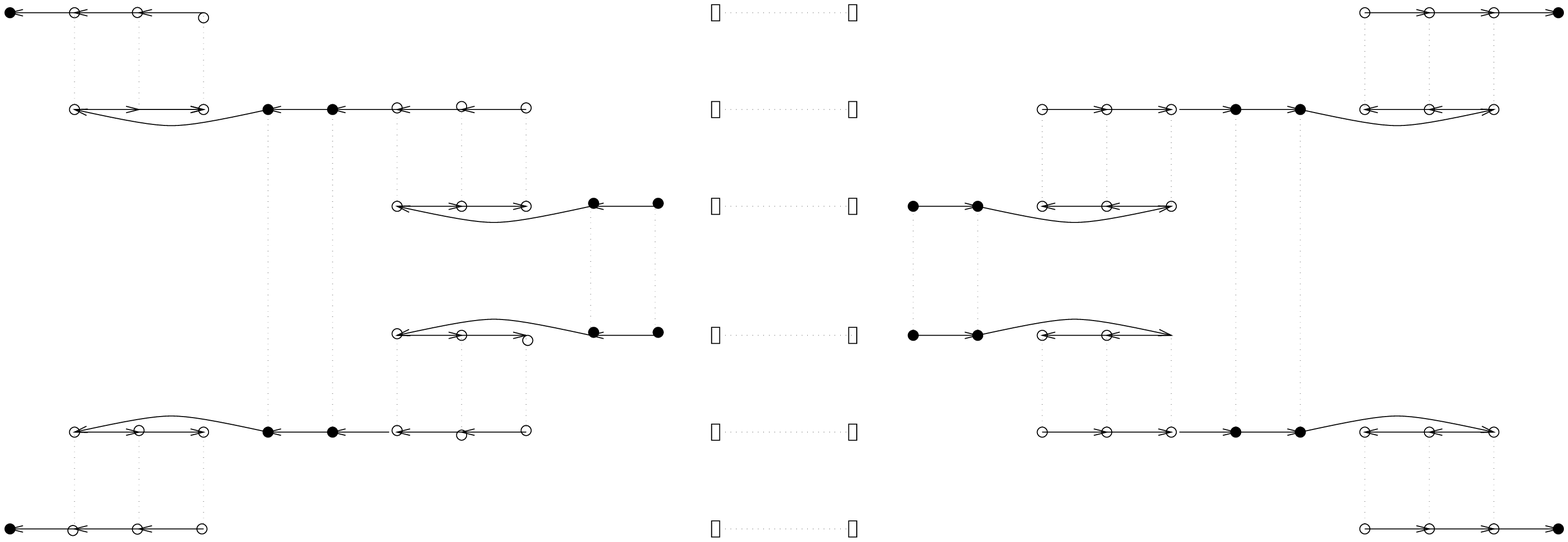} 
\end{figure}
\end{example}

What we have is called a bunch of chains ~\cite{mp}, which 
 codes a matrix problem. 
The number of matrices can be infinite, but this does not disturb the 
general theory.  

Namely, the matrix problem is the following:
\begin{enumerate}
\item Each triple $(L,a,i)$, ($L$ is a component of $\tilde{X}$, $a\in L$,
$i\ge 0$ the integer number) corresponds to some matrix $M(L,a,i)$. 
These matrices are divided into horizontal blocks, numbered by the points 
of $E_{L,a}(i)$.
Since $F_{L,a}$ consists only of one element, we do not have a vertical division in this case. Ideed, some blocks may have zero size, i.e., be empty.
\item Blocks, corresponding to the conjugated points from $\bf E$, 
        have the equal number
       of rows; points, corresponding to the conjugated points from $\bf F$,
 have the equal number of columns
\item We have the partial order on the set of points. Let us say that the horisontal blocks are supplied with some ``weights'', and the weight of  one block is bigger than the weight of another block if the point, corresponding to the first block, is bigger than the point corresponding to the second one.

\item We can do the following transformation with our matrices:

\begin{enumerate}
\item 
Simultaneous:  elementary transformations with the columns of matrices $M(L,a,i)$ and $M(L',a,i)$.
\item Simultaneous: any elementary transformations inside of the conjugated blocks.
\item Independent: add a scalar multiple of any row from the block with lower 
weight to any row of a block of the higher weight.
\end{enumerate}
\end{enumerate}

In our case there are some additional restrictions on our matrices.
\begin{enumerate}
\item All big matrices are square and nondegenerated.
\item If one of the conjugated blocks is nonempty, that the other one 
      is nonempty, too.
\end{enumerate}

There are two types of indecomposable objects: bands and strings.

\begin{enumerate}
\item Band data $\kB(w,m,\lambda)$ is given by two discrete parameters:
by word $w$, natural number  $m$, and one continuous parameter 
$\lambda \in k^{*}$.
A word $w$ is just a sequence of points of ${\bf E} \cup {\bf F}$ $w_{1} - w_{2} \sim w_{3} - w_{4} \sim \dots - w_{N}$, connected  by the symbols of two types,
$-$ and $\sim$. The symbol $\sim$ should stay between conjugated points,
 $-$  only between
a point of the type $E_{L,a}(*,i)$ and a point $F_{L,a}(i)$. If one link was 
$-$, then the next one should be $\sim$ and vice versa.  In a band data a word 
$w$ should be closed: $x_{N} \sim x_{1}$. It means that it can be written as a cycle. We require that $w$ is not a power of some other word. 

\item A string data $\kT(w)$ depends only on some  full word $w$. Full means
that  $w$ contains each point  $w_{i}$ together with its conjugate. In our
situation a word $w$ could be infinite. We require, however, that each
point $w_{i}$ appears only a finite number of times. 
\end{enumerate}

Let us briefly recall the algorithm, giving a concrete description of matrices, corresponding to the band and string data.

\begin{enumerate}
\item Let  a band data be $\kB(w,m,\lambda)$.
We count the entrance of each class of conjugated points. Let point $w_{i}$ 
occured $k_{i}$ times. The block,
corresponding to $w_{i}$,  should be divided into $k_{i}$ strips. We have a division of a big matrix $M(L,a,i)$ into the smaller blocks. Let us look
now at the subwords $w_{i} - w_{i+1}$. Suppose we have the $k$-th appearance of the class $[w_{i}]$ and $l$-th appearance of the class $[w_{i+1}]$. One of 
the points $w_{i}$,$w_{i+1}$ belongs to ${\bf E}$, the other to ${\bf F}$.
If $w_{i+1} \ne  w_{N}$, then we put in the entry with the coordinates 
$(k,l)$ (with respect to the subpartition of  $w_{i} \times  w_{i+1}$ 
submatrix) the identity matrix $I_{m}$ (here our second discrete parameter 
appears). If $w_{i+1} =  w_{N}$, then we put, on the corresponding place, the 
Jordan block $J_{m}(\lambda)$. All other entries are zero. 

\item Let  a string data be $\kT(w)$. The algorithm is basically the same as in the case of bands. The only difference is that we have to put not 
the $I_m$ or $J_{m}(\lambda)$, but just $1\times 1$ matrix $1$.

\end{enumerate}

It is clear that the closed word is necessarily finite. A word, defining a
string data, could be infinite.

Consider an example. We  write down the data, corresponding to 
the complexes on the nodal rational curve, which were considered in the 
previous chapter.

\begin{example} Let $w = $
$\stackrel{1}E_{(X,a)}(0,-1\omega_{0}) - F_{(X,a)}(0) 
\stackrel{1}\sim F_{(X,b)}(0) - E_{(X,b)}(0, 1\omega_{0})
\stackrel{1}\sim E_{(X,a)}(0,1\omega_{0}) - F_{(X,a)}(0) 
\stackrel{2}\sim F_{(X,b)}(0) - E_{(X,b)}(0,-2\omega_{1}) 
\sim E_{(X,b)}(1,2\omega_{-1}) - F_{(X,b)}(1) \stackrel{1}\sim F_{(X,a)}(1) - 
E_{(X,a)}(1,\omega_{-1})$ $
\stackrel{1}\sim E_{(X,a)}(0, -\omega_{1}) - F_{(X,a)}(0) \stackrel{3}\sim F_{(X,b)}(0)
- E_{(X,b)}(0, -4\omega_{1}) \stackrel{1}\sim E_{(X,b)}(1,4\omega_{-1}) - F_{(X,b)}(1) \stackrel{2}\sim
F_{(X,a)}(1) - E_{(X,a)}(1,0\omega_{0}) \stackrel{1}\sim E_{(X,b)}(1,0\omega_{0}) - $ $
F_{(X,b)}(1) \stackrel{3}\sim F_{(X,a)}(1) - E_{(X,a)}(1,-2\omega_{-1}) 
\stackrel{1}\sim E_{(X,a)}(2,2\omega_{-1}) -
 F_{(X,a)}(2) \stackrel{1}\sim F_{(X,b)}(2)  - E_{(X,b)}(2, -2\omega_{0}) 
\stackrel{1}\sim
E_{X,a}(2, -2\omega_{0}) - F_{(X,a)}(2) \stackrel{2}\sim F_{(X,b)}(2) - $ 
$ E_{(X,b)}(2,
1\omega_{-1}) \stackrel{1}\sim E_{(X,b)}(1, -1\omega_{1}) - F_{(X,b)}(1) 
\stackrel{4}\sim F_{(X,a)}(1)
- E_{(X,a)}(1,1\omega_{-1}) \stackrel{2}\sim E_{(X,a)}(0,-1\omega_{-1}) - F_{(X,a)}(0) \stackrel{4}\sim
F_{(X,b)}(0) - E_{(X,b)}(0,-\omega_{0}) \stackrel{1}\sim
$.
Then $\kB(w,n,\lambda)$ is just a complex from Example 3.4 from the 
previous chapter.
\end{example}

\begin{theorem}(See ~\cite{mp}.)
\begin{enumerate}
\item
All representations $\kB(\omega, m, \lambda)$, $\kS(\omega)$ are indecomposable.
Each  indecomposable representation is isomorphic,
either to some band representation $\kB(\omega, m, \lambda)$ or
to some string representation $\kS(\omega)$.

\item The only isomorphisms between these objects are
\begin{enumerate}
\item $\kS(\omega) \cong \kS(\omega^{-1})$, where 
$\omega = a_0 r_1 a_1 \dots r_m a_m$ and \
$\omega^{-1} = a_m r_m a_{m-1} \dots r_1 a_0$ is the inverse word.
\item $\kB(\omega, m, \lambda)= \kB(\omega', m, \lambda')$, where $\omega'$ 
is a cyclic permutation of
 $\omega$  and $\lambda'$ depending on the signum of the permutation
being either
$\lambda$ or $\lambda^{-1}$.
\end{enumerate}
\end{enumerate}
\end{theorem}

We want now to illustrate the convenience of our description of the complexes
in the derived category $D^{-}(\Coh_{X})$. Let $\kF_{.}$ and $\kG_{.}$ be two
objects, given by triples $(\tilde{\kF}_{.},\kM, i)$ and $(\tilde{\kG}_{.},
\kN_{.},j)$. We can ask  which triple corresponds to the tensor
product of complexes $\kF_{.} \tens_{L} \kG_{.}$. As one can easily see,
it should be 
$
(\tilde{\kF}_{.}\otimes \tilde{\kG}_{.}, \kM_{.}\otimes \kN_{.},
i\otimes j).
$

By the K\"unneth formula we have a functorial isomorphism 
(since the homological dimension is $0$):
$$
 \bigoplus\limits_{k+l=n} (H_{k}(\kM_{.})\otimes (H_{l}(\kN_{.}) 
\xrightarrow{\oplus (H_{k}(i)\otimes H_{l}(j))} 
  H_{n}(\kM_{.}\otimes \kN_{.}).
$$

This means that we can compute the matrices, corresponding to the tensor
product of complexes.

\section{\bf Description of the coherent sheaves, vector bundles, torsion-free sheaves, mixed sheaves and skyscraper sheaves}

Now we want to show what corresponds to coherent sheaves. Let complex 
$\kF_{.}$ be given by a triple $(\tilde{\kF}_{.}, \kM_{.}, i)$. 
We have to write the conditions \ $H_{i}(\kF_{.})=0, (i \ge 1)$ in the language of matrices. Recall that we have the following diagram:

\begin{tabular}{p{1.4cm}c}
&
\xymatrix
{
{\mathcal J}\tilde{\mathcal F}_{.} \ar[r] &
\tilde{\mathcal F}_{.} \ar[r] &
\tilde{\mathcal F}_{.} \tens_{\tilde \kO} \tilde{\kO}/{\mathcal J} \ar[r]&
{\mathcal J}\tilde{\mathcal F}_{.}[-1] \\
{\mathcal J}\tilde{\mathcal F}_{.} \ar[r] \ar[u]^{id} &
{\mathcal F}_{.} \ar[r] \ar[u]_{\Phi} &
{\mathcal M}_{.} \ar[r] \ar[u]_{i} &
{\mathcal J}\tilde{\mathcal F}_{.}[-1]. \ar[u]^{id}\\
}
\end{tabular}

Write the long exact sequence of homologies, associated with this morphism
of triangles:
$$
\xymatrix
{
0 & 
H_{0}(\tilde{\kF_{.}}/\kJ \tilde{\kF_{.}}) \ar[l] &
H_{0}(\tilde{\kF_{.}}) \ar[l] & 
H_{0}(\kJ \tilde{\kF_{.}}) \ar[l] &
H_{1}(\tilde{\kF_{.}}/\kJ \tilde{\kF_{.}}) \ar[l] &
\dots \ar[l] \\
0 &
 H_{0}(\kM_{.}) \ar[l] \ar[u]^{H_{0}(i)} &
 H_{0}(\kF_{.}) \ar[l] \ar[u]  &
 H_{0}(\kJ \tilde{\kF_{.}}) \ar[l] \ar[u]^{id} &
 H_{1}(\kM_{.}) \ar[l] \ar[u]^{H_1(i)} &
\dots \ar[l] \\
}
$$

\noindent
We get
$H_{i}(\kF_{.})=0, (i \ge 1)$ is equivalent to 
\begin{enumerate}
\item $H_{1}(\kM_{.})\lar H_{1}(\tilde{\kF}_{.}/\kJ\tilde{\kF}_{.} ) 
\lar H_{0}(\kJ\tilde{\kF}_{.})$
being a monomorphism.
\item 
$H_{k+1}(\kM_{.})\lar H_{k+1}(\tilde{\kF}_{.}/\kJ\tilde{\kF}_{.} ) 
\lar H_{k}(\kJ\tilde{\kF}_{.})$
being an isomorphism. 
\end{enumerate}

Let us give a combinatorical interpretation of these conditions (for both bands and strings). Let $w$ be a parameter either of $\kB(w,m,\lambda)$ or of
$\kS(\omega)$. These conditions imply:

\begin{enumerate}
\item  A word  $w$  does not contain any 
$E_{L,a}(k+1, n\omega_{0}), k\ge 0$; 
\item 
a word $w$ does not contain any  
$E_{L,a}(k+1, n\omega_{-1}) (n\ge 1, k\ge 1)$, $E_{L,a}(k+1, -m\omega_{1}),
(m\ge 2, k\ge 1)$;
\item
a word  $w$ does not contain any subword of type $ E_{L,a}(k+1, \omega_{1}) - F_{L,a}(k+1) \sim F_{L',a}(k+1)
- E_{L',a}(k+1, \omega_{1})$.
\end{enumerate}

\noindent
It gives us the description of coherent sheaves:

\begin{enumerate}
\item All the bands $\kB(w,m,\lambda)$ such that the word $w$ doesn't contain
 points $F_{L,a}(i)$ with $i \ge 2$. 
     
\item Strings  $\kT(w)$ with the following properties:
\begin{enumerate}
\item There are no points $E_{L,a}(n\omega_{0},i)$ with $i \ge 2$,
       $E_{L,a}(k+1, n\omega_{-1})$, $E_{L,a}(k+1, -m\omega_{1})$ with 
    $n,m \ge 2$.
\item $w$ does not contain any subword of type 
       $ E_{L,a}(k+1, \omega_{1}) - F_{L,a}(k+1) \sim F_{L',a}(k+1)
       - E_{L',a}(k+1, \omega_{1})$.
\item For each $i \ge 0 $ the points of the type $E_{L,a}(*,i)$ appeared the 
same number of times as $F_{L,a}(i)$ ( it is a condition for the matrices
$M(L,a,i)$ to be square and nondegenerated).
\end{enumerate}
\end{enumerate}

In a similar way we can describe the  {\it bounded} derived category 
$D^{b}(Coh_{X})$: the conditions $H_{i}(\kF_{.}), i \gg 0$ can be described
in a similar way. 

In particular we get a description of 

\begin{enumerate}
\item Vector bundles (we get just the matrix problem from the work of
 Drozd and \ Greuel): bands $\kB(w,m,\lambda)$
with $w$ not containing $F_{L,a}(i), (i \ge 1)$ (in case of a curve of 
arithmetic genus $1$). 

\item Skyscraper sheaves: bands and strings, defining a coherent sheaf and not containing $E_{L,a}(n\omega_{0},i), i\ge 0$ (it follows from the observation
 that $\kF$ and $\kF\tens_{\kO}\tilde{\kO}$ have the same support).
\end{enumerate} 

We see that for a coherent sheaf $\kF$ we have either 
${\mathcal Tor}_{\kO}^{i}(\kF, \kO/\kJ) = 0$ for $i > 1$ or it is nonzero
for all $i \ge 2$. 
As a corollary we obtain that the homological dimension
of an object of $Coh_{X}$ is either $0$ or $1$ or $\infty$ (which coincides 
with the result of the Auslander-Buchsbaum formula).

We are going to do the last step in our classification: we shall describe among
 all the coherent sheaves, the torsion-free sheaves. Those which are not vector bundles have infinite homological dimension, and hence we should look for them among strings.  Let $\kF$ be a coherent
sheaf on $X$. It is torsion-free if and only if all its localizations
$\kF_{x}$ are torsion-free $\kO_{X,x}$-modules. At the regular point this 
condition is obvious. But we go further, 
$\kF_{x}$ is a torsion-free $\kO_{X,x}$-module if and only if its completion
$\hat{\kF_{x}}$ is a torsion-free $\hat{\kO}_{X,x}$-module. But in our case, if
$x$ is singular, then $\hat{\kO}_{X,x} = k[[x,y]]/(xy)$. The indecomposable 
torsion-free
modules are known in this case, they are either $k[[x]]$ or $k[[y]]$
or the regular module  $k[[x,y]]/(xy)$ ~\cite{Bass}. 

Now, let us mention  that in the same way we have dealt with curves, we can
deal with the local ring $k[[x,y]]/(xy)$. Namely, we  consider its 
normalization $k[[x]]\times k[[y]]$, conductor $J=(x,y)$ and just repeat
the construction of the category of triples. 
As a result, we get the following matrix problem 
(in the notation of Bondarenko):

\begin{figure}[ht]
\hspace{4cm}
\includegraphics[height=3.5cm,width=4cm,angle=-90]{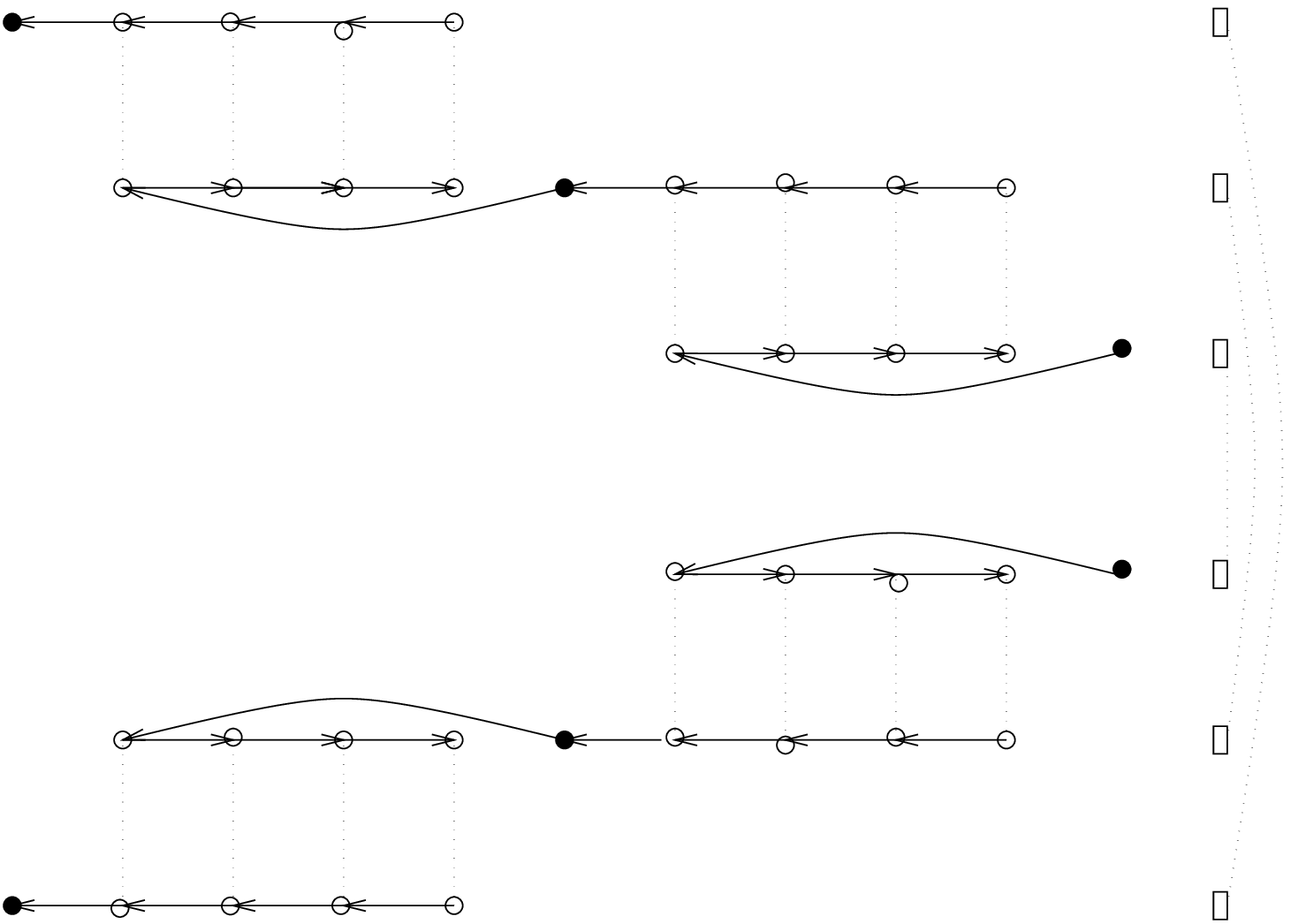}
\end{figure}

Let $x\in X$ be singular. Consider the functor $\Coh_{X} \lar (\kO_{X,x}-\mod) 
\lar (k[[x,y]]/(xy)-\mod)$ (composition of the localization and completion).
This functor is exact and so induces the functor between the derived categories
$D^{-}(\Coh_{X}) \lar D^{-}(k[[x,y]]/(xy)-\mod)$. 
What does it look like on triples? Obviously, $(\tilde{\kF}_{.},\kM_{.},i)$
is mapped to $\widehat{(\tilde{\kF}_{.})}_{x}, \widehat{(\kM_{.})}_{x}, i_{x})$.
So, the image of the triple is described by the same matrices! But surely
there is one important difference: blocks, corresponding to  vector
bundles are united, and there are no links between them anymore.
But we know how the modules $k[[x]]$, $k[[y]]$ and $k[[x,y]]/(x,y)$ are given 
in the language of triples. Let us denote $T_{nx} : k[[x]] \stackrel{x^n}\lar
k[[x]] $,$T_{ny} : k[[y]] \stackrel{y^n}\lar k[[y]]$. 
Then $k[[x]]$, for example, is given by 
normalization: 
$$k[[x]] \oplus (\bigoplus\limits_{i=0}^{\infty} T_{iy}[-i] 
\oplus T_{i1 x}[-i-1])$$
and matrices: 
\begin{figure}[ht]
\hspace{4.3cm}
\includegraphics[height=3.5cm,width=2.4cm,angle=-90]{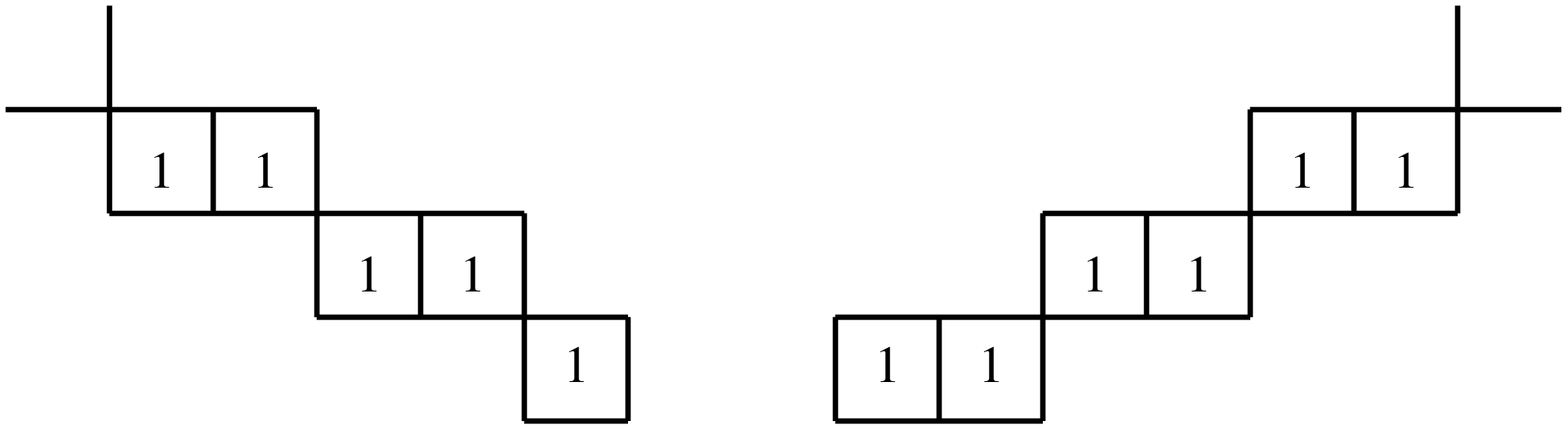}
\end{figure}

Hence we can deduce the answer:
torsion-free sheaves, which are not vector bundles, are strings
 $\kT(w)$, where  $w$ does not contain any  $E_{L,a}(i, n\omega_{-1})$, 
\ $E_{L,a}(i, -m\omega_{1})$
$(i\ge 2, m,n \ge 2)$. Moreover, each $E_{L,a}(i, \omega_{-1})$, $E_{L,a}(i, -\omega_{1})$ 
can occur in word $w$ at most one time.   

\section{\bf The construction of Polishchuk}

In a recent paper A.~Polishchuk ~\cite{Polishchuk} showed a connection between
the structure of the derived category of the coherent sheaves 
$D^{b}(\Coh_{X})$,
where $X$ is a projective curve of arithmetical genus $1$ with nodal singularities and the trihonometric solutions of the classical Yang-Baxter equation. The 
main role in this construction is played by the so called spherical objects:

\begin{definition}(See ~\cite{Thomas}.)
Let  ${\mathcal D}$ be a triangulated category over a field  $k$, such that all 
spaces 
$\Hom(X,Y)$ are finite-dimensional.
$\Hom^{i}(M,N) := \Hom(M, N[i])$. An object $M$ is called $n$-spherical, if
\begin{enumerate} 
\item
$\Hom^{i}(M,M)=0 \ \forall i \ne 0,n$, 
$\Hom^{0}(M,M)\cong \Hom^{n}(M,M)\cong k$.
\item $\forall F \in Ob({\mathcal D})$ the composition map 
$\Hom^{i}(M,F)\times \Hom^{n-i}(F,M) \lar \Hom^{n}(M,M) \cong k$ is 
nondegenerated.
\end{enumerate} 
\end{definition}

Let  $X$ be a projective curve as above.
${\mathcal D}= D^{b}(Coh_{X})$. It  
 was shown in ~\cite{Polishchuk} that all simple 
vector bundles are  1-spherical.
Indeed, the first condition follows from the theorem of Riemann-Roch for singular curves. Namely, let $\kE$ be a simple vector bundle. This implies that
$\Hom(\kE, \kE) \cong k$. We have to show that  
$\Ext(\kE, \kE) \cong k$. From the local-global spectral sequence 
follows that 
$$
\Ext(\kE, \kE) \cong H^{1}({\mathcal Hom}(\kE, \kE)) \cong H^{1}(\kE^{\vee} \tens \kE).$$
$\kE^{\vee} \tens \kE$ is a vector bundle. By the theorem of Riemann-Roch 
$\chi(\kE^{\vee} \tens \kE) = deg(\kE^{\vee} \tens \kE) + 1 - p_{a}(X) = 0$
holds. Hence
 $\dim_{k}(\Ext(\kE, \kE))=\dim_{k}(\Hom(\kE, \kE))=1$.

The second condition follows from the Serre duality for the derived categories
(see for example ~\cite{Thomas}).
One can also show that the skyscraper sheaves
 $\kO_{x}/m_{x}$ ($x$ ist regular) are also 1-spherical.

\begin{conjecture}(See ~\cite{Polishchuk}.) What are the  1-spherical 
objects in this case? Which orbits have
the set of spherical objects under the action of the group of autoeqivalences of the derived category $Aut({\mathcal D})$?
\end{conjecture}

\vspace{1cm}
\begin{tabular}{lp{1.5cm}l}
Igor Burban & & Yuriy Drozd\\
Universit\"at Kaiserslautern & & Kyiv University \\
burban@mathematik.uni-kl.de & & yuriy@drozd.org
\end{tabular}


\begin{thebibliography}{99}

\bibitem{At} M.~Atiyah, {\it Vector bundles over an elliptic curve}, Proc. London
Math. Soc., {\bf 7}(1957), 414-452.

\bibitem{Bass} H.~Bass, 
{\it On the ubiquity of Gorenstein rings}, Math. Zeitsch.,
{\bf 82} (1963), 8-27.


\bibitem{mp} V.~V.~Bondarenko, {\it Representations of bundles of semi-chains and their applications}, St. Petersburg Math. J., {\bf 3} (1992), 973-996


\bibitem{mp1} V.~V.~Bondarenko, {\it Bundles of semi-chains and their representations}, preprint of the Kiev Institute of mathematics, 1988.

\bibitem{Roiter2}
V.~V.~Bondarenko, L.~A.~Nazarova, A.~V.~Roiter, V.~V.~Sergijchuck
{\it Applications of the modules over a diad to the classification of
finite $p$-groups, having an abelian subgroup of index $p$},
Zapiski Nauchn. Seminara LOMI,
{\bf 28} (1972) 69-92.


\bibitem{vb} Yu.~A.~Drozd, G.-M.~Greuel, {\it On the classification of the 
vector bundles on projective curves}, Max-Plank-Institut f\"ur Mathematik, 
Preprint Series
1999 (130).

\bibitem{bimproblems} Yu.~A.~Drozd, 
{\it Matrix problems and categories of matrices}, Zapiski Nauchn. Seminara LOMI,
{\bf 28} (1972) 144-153.

\bibitem{Gelfand} I.~M.~Gelfand, {\it Cohomology of the infinite dimensional
Lie algebras; some questions of the integral geometry}, International congress
of mathematics, Nice, 1970.

\bibitem{GP} I.~M.~Gelfand, V.~A.~Ponomarev, {\it Indecomposable representations
of the Lorenz group}, Uspehi Mat. Nauk, 1968, {\bf 140}, 3-60.

\bibitem{Groth} A.~Grothendieck, {\it Sur la classification des fibres holomorphes
sur la sph\`ere de Riemann}, Amer. J. Math., {\bf 79} (1956), 121-138.

\bibitem{Hartshorne}
Hartshorne,~R.: Algebraic Geometry, Springer-Verlag (1977). 

\bibitem{Roiter1}
L.~A.~Nazarova, A.~V.~Roiter, {\it Finitely generated modules over diad 
of two discrete valuation rings}, Izv. Akad. Nauk USSR, ser. mat. {\bf 33}, 
(1969), 65-89.

\bibitem{Roiter3}
L.~A.~Nazarova, A.~V.~Roiter, {\it About one problem of I.~M.~Gelfand},
Functional analysis and its applications, vol. {\bf 7}, {\bf 4}, 1973, 54-69.


\bibitem{Polishchuk} 
A.~Polishchuk, {\it Classical Yang-Baxter equation and the $A_{\infty}$-constraint}, arXiv:math.AG/0008156. 

\bibitem{Thomas} P.~Seidel, R.~P.~Thomas,
{\it Braid group actions on derived categories of coherent sheaves},
arXiv:math.AG/0001043.  

\end{thebibliography}
\end{document}